# Logic: a Misleading Concept.
# A Contradiction Study toward Agent's Logic Ontology

Feng Liu
Department of Management Science and Engineering
Shaanxi Economics and Trade Institute (South Campus)
South Cuihua Road, Xi'an, Shaanxi, 710061, P. R. China

Florentin Smarandache
Department of Mathematics
University of New Mexico, Gallup, NM 87301, USA

**Abstract:** The paper presents a fresh new comprehensive ideology on Neutrosophic Logic based on contradiction study in a broad sense: general critics on conventional logic by examining the essence of logic, fresh insights on logic definition based on Chinese philosophical survey, and a novel and genetic logic model as the elementary cell against Von Neumann oriented ones based on this novel definition. As for the logic definition, the paper illustrates that logic is rather a tradeoff between different factors than truth and false abstraction. It is stressed that the kernel of any intelligent system is exactly a contradiction model. The paper aims to solve the chaos of logic and exhibit the potential power of neutrosophy: a new branch of scientific philosophy.

1. **Background**

Although it is commonly believed that intelligence is a social activity, and it is therefore represented in multiagent forms, but its kernel, the logic of agents, remains controversial with its static, monolateral or homogeneous forms.

This reflects in their behaviors as: our agents appear social in outer forms but autarchic in nature, for this kind of multiagent system can never deal with controversies, critics, conflicts or something with flexibility. Our multiagent system has become a sort of software engineering or system engineering of fresh forms, failing to implement our presumed social intelligence.

In the long-term exploration, one realizes that the problem takes its root in the misleading definition of logic. Even the simplest logic such as "The earth turns around the sun" and "I'll visit him if it doesn't rain and he is in" can lead to ambiguous or contradictory actions of agent (Liu, [2]). Limited to the length, I'll present in this paper only a framework to launch our discussion, as follows:
- Fact: a belief rather than truth
- Logic: dependent of situations, not absolute
- Logic is negating itself
- Logic is only one perspective of learning, not an independent entity
- As a part of learning, logic is dynamic
- As a part of learning, logic is multilateral
- Logic is always partial
- Illusion and creativity

Many scientists argue about the need to model human intelligence in the general level. The argument lies in our vague understanding of intelligent system (Liu, [1]). Intelligent system should be, in our opinion, **a tradeoff machine in order to adapt to its environment**. Then a specific model becomes such a tradeoff between ideal philosophic model and practical system model, in the hierarchy from philosophic layer down to a specific application or situation constraint implementation. I'll show this philosophy in step ward way, as follows.

## 2. Neutrosophy

Neutrosophy is a new branch of philosophy that studies the origin, nature, and scope of neutralities, as well as their interactions with different ideational spectra.
It is the base of *neutrosophic logic*, a multiple value logic that generalizes the fuzzy logic and deals with paradoxes, contradictions, antitheses, antinomies.

**Characteristics** of this mode of thinking:
- proposes new philosophical theses, principles, laws, methods, formulas, movements;
- reveals that world is full of indeterminacy;
- interprets the uninterpretable;
- regards, from many different angles, old concepts, systems: showing that an idea, which is true in a given referential system, may be false in another one, and vice versa;
- attempts to make peace in the war of ideas, and to make war in the peaceful ideas;
- measures the stability of unstable systems, and instability of stable systems.

Let's note by <A> an idea, or proposition, theory, event, concept, entity, by <Non-A> what is not <A>, and by <Anti-A> the opposite of <A>. Also, <Neut-A> means what is neither <A> nor <Anti-A>, i.e. neutrality in between the two extremes. And <A'> a version of <A>.
  <Non-A> is different from <Anti-A>.

**Main Principle:**
Between an idea <A> and its opposite <Anti-A>, there is a continuum-power spectrum of neutralities <Neut-A>.

**Fundamental Thesis of Neutrosophy:**
Any idea <A> is T% true, I% indeterminate, and F% false, where T, I, F $\subset\ ]^-0, 1^+[$.
Here $]^-0, 1^+[$ is a non-standard unit interval, with $^-0=\{0-\varepsilon, \varepsilon$ is a positive infinitesimal number$\}$ and $1^+=\{1+\varepsilon, \varepsilon$ is a positive infinitesimal number$\}$.

**Main Laws of Neutrosophy:**
Let <α> be an attribute, and (T, I, F) $\subset\ ]^-0, 1^+[^3$. Then:
- There is a proposition <P> and a referential system {R}, such that <P> is T% <α>, I% indeterminate or <Neut-α>, and F% <Anti-α>.
- For any proposition <P>, there is a referential system {R}, such that <P> is T% <α>, I% indeterminate or <Neut-α>, and F% <Anti-α>.
- <α> is at some degree <Anti-α>, while <Anti-α> is at some degree <α>.

## 3. Fact: a Belief rather than Truth

We start with an ancient problem based on the following contradiction:
- The sun turns around the earth.
- The earth turns around the sun.

Of cause nearly everyone of us would answer: the later is absolute right. Note that this is merely a belief, because in Copernicus's age the majority believed in the former. Has anyone proved nowadays whether the former is incorrect? If yes, he must have assumed that the sun is relatively fixed. Unfortunately this is also his belief, because none of us has ever proved the absoluteness of his consciousness: when we see something, is it really something or just we believe that there is something (we really touch something or we really believe it is something we touched)? Or more specifically, is it an object or just we hold long this same belief? Do we really exist as in form we see or just we believe so? I have to introduce a heard experiment to show this point.

> A blindfold person is told to be experimented with an iron burnt hot. And through a chronic preparation before him, the iron is burnt fervid, and he is told that the iron is gradually moved closer and closer to him.
> "Yes, I am feeling hotter and hotter, …… really hot, extremely, ……"

The gradual process goes on and on, until suddenly, he is instructed to have his skin burnt.
"Oh……", his skin really burnt.
When he opened his eyes, there is nothing but the scorch in him—there is no fire nor iron, but merely his **imagination—it is strong enough to cause the effect.**

I experienced another experiment in which four of us were pointing to a carefully set small wooden stool while rotating around it. According our mutual will, the stool turned itself in the same direction we turned!

Another fact (shown in a qigong journal quite a number of years ago, the following is based on our memory) shows the same thing:

There is a qigong (commonly believed as some mental or physical exercise in order to gather the "energy" from nature, qi (there are such a kind of substance in Chinese medicine which is unseen but really affects our body), or the concentrative power of will to maintain health from disease, it is not a feasible way to us) expert in China who, through chronic practice, can "brake" a steel saw blade with nothing but his will, and he had been succeeding in it nearly every time, even in many qigong reports.

Once he re-showed the same talent to the huge audience with great curiosity. He ordered: "break", but unexpectedly, the blade remain exact the same as before, and the following tries turned out to be the same failures. The atmosphere became extremely unfavorable.

Fortunately however, the chairman of the qigong report is experienced, and asked the audience to cooperate: more you are confident, more successful the experiment.

Magically, the expert broke the blade with a single command.

Conclusion:

We have to confess from the three experiments that **fact is really our belief**—if a single belief is not powerful enough to convince us, the mutual belief, especially of all the human beings, is definitely strong enough to illude ourselves. While this cause-effect goes on and on, we are unconsciously trapped in the inextricable web of deceit designed by ourselves.

Only wise man can see through this kind of deceit, e.g., Master Huineng in Chinese Tang dynasty when he saw an argument about a pennant aflutter: whether the wind was moving or the pennant.

**"Neither. What actually moved were your own minds."**(Liu [4], see also Yan Kuanhu [2])

Everyone can become wise when understands this cause-effect, which is the basic point of Buddhism (Chin Kung [1])

4. **Logic: Dependent of Situations, not Absolute**

Take the logic 1+1=2 for example. Is it correct? Consider
black+white=?, explosive+fire=?, warm+cold=?, theory+practice=?, and yin+yang = ?

Does the idiom "Blind People Touching an Elephant" really refer to blind men and elephant?

- Blind People Touching an Elephant, a story from the Mahapra Janaparamita Sutra:
  The story shows that **the same elephant can be interpreted as such different things** as turnip, dustpan, pestle, bed, jar and rope **by different blind people** who touch it in turn. The first one touches the tusk, the second the ear, the third the foot, the fourth the back, the fifth the belly, and the last the tail.
  **Based on their different beliefs, the same elephant conveys diverse logics**.

Conclusion:

Logic is more a kind of mental behavior than an objective understanding, i.e., it is more a belief rather than truth; this belief is based on "facts" which are also beliefs.

This belief is subject to dynamic changes with situations, and more general belief (general understanding) relative to more general situations could be too flexible to grasp (e.g., Dao=yin+yang), therefore **logic suggests varying explanations based on particularity of situations.**

5. **Logic is Negating Itself**

Logic comes as mental reflection and leads to new reflection. So, **it is not the problem of logic (validity) but the ways we reflect it**, otherwise it would become the Chinese room experiment.

- J. R. Searle shows this "Chinese room problem" in his paper "Minds, brains and programs":

We set an Englishman which does not know Chinese, in a closed room, with many symbols of the Chinese language, and a book of instructions in English of how to manipulate the symbols when a set of symbols (instructions) is given. So, Chinese scientists will give him instructions in Chinese, and the **Englishman will manipulate symbols in Chinese, and he will give a correct answer in Chinese. But he is not conscious of what he did.** We suppose that a machine behaves in a similar way: it might give correct answers, but it is not conscious of what it is doing.

Another argument on validity of logic is based on a Chinese idiom: Cutting a Mark on the Boatside to Retrieve a Sword (Young):

> Once, a man of the State of Chu (ancient China) took a boat to cross a river. It so happened that his sword slipped off and fell into the water. Immediately he cut a mark on the side of the boat and assured himself: "This is where I have dropped my sword."
>
> By and by the boat came to the destination and stopped. The man plunged into the stream at the point indicated by the incised mark trying to retrieve the lost sword.
>
> The boat has moved on, but not the sword. To recover his sword this way—the man is indeed muddle-headed .

This prompts us to doubt whether logic is always applicable to other circumstances as we know situation is subjective to constant change. It can be successful in closed systems where every state is well defined, but how about open ones?

Daodejing (Wang Bi, Guo Xiang) begins with: "Dao, daoable, but not the normal dao." Referring to the natural law, we can say it is dao, but it doesn't mean what we say. Whenever we mention it, it is beyond the original sense.

In Daodejing the creator of everything is defined as dao: like a mother that bears things with shape and form. But what/who is dao? It is just unimaginable, because whenever we imagine it, our imagination can never be it (we can never completely describe it: more we describe it, more wrong we are). Daoism illustrates the origin of everything as such a form that doesn't show in any form we can perceive. Whatever we can perceive is merely the created forms, rather than its genuine nature, as if we know people by their outer looks rather than by their inner intentions. We are too far from understanding the nature.

**Daodejing suggests that logic in the most original extremity is shapeless in nature:** it is unbodied, invisible, inexpressible, or even intangible.

- We frequently have such a feeling in learning English as a foreign language that there is no fixed meaning but an intangible impression or feeling to a word: the meaning varies with situations, contexts or even ages so that we can never assure our comprehension. In fact, it is due to the unbridled usage in logic made by people of different ages and districts—there is always creativity implied in the word so that we can rely on nothing more than our own creativity.
- Wherever there is logic, there is also the corresponding comprehensive understanding to logic: also based on our creativity relative to the different existences of human beings. Due to the diversity of comprehension, interpretation and creativity to the same logic, there are varying versions of perception (or conception) to the same logic.

Whenever we say "it is" by logic, we are subjective—how can we assure our objectiveness? We may have developed it from some strictly limited domain under constrained conditions, or it is merely a haphazard, because we can only observe limited cases in our limited lives. How can we convey the logic to those with different backgrounds, even with the slightest difference? How can we then assure the determinacy of its truthness? Is that the reason that majority of people hold it (e.g., Darwin's evolution theory and the functional difference between left and right brain—both are controversial in fact)? In fact none of the logics can be proved, even of we exist or not (Gershenson).

There is no truth and false actually: there is because the outcome has to meet someone's desire—they are merely the attributes of a tradeoff. One false dead can be true in another perspective, e.g., eating much is good, because of the excellent taste and nourishment, but it is also bad when he gets weighted. Neutrosophy (Smarandache) shows that a true proposition to one referential system can be false to another.

Conclusion:

Validity of logic depends on the way we reflect it, not logic itself. Logic never proves itself.

**Logic is a matter of tradeoff (balance) between contradictory factors. There seems no absolute correctness or falseness independent of environment.**

There is dao, but not the kind we mentioned, **accordingly, there is logic but not what we specified**.

## 6. Logic is only One Perspective of Learning, not an Independent Entity

Logic comes from perception and leads to new perception. It is shown that **human understands the world through the interaction of the inter-contradictory and inter-complementary two kinds of knowledge: perceptual knowledge and rational knowledge——they can't be split apart.**

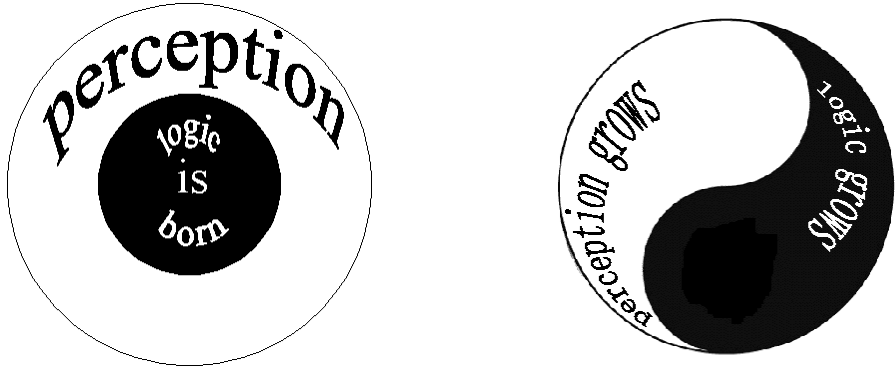

   Logic is created through perception in which void, intangible feelings or impressions have been nurtured, brought up and developed into a mental model.
   The logic born is nothing more than a subjective hypothesis at primitive stage—it is not within the sense of truthness and falseness.
   However, the terms (symbols) and syntax (rules) is only understood by perception through practice. How can we imagine a bookworm who well reads books but has no experience?
   Where are truth and false born? There are no such beliefs at the first stage of practice in fact (Daodejing):

> When beauty is abstracted (Peter A. Merel)
> Then ugliness has been implied;
> When good is abstracted
> Then evil has been implied.
>
> So alive and dead are abstracted from nature,
> Difficult and easy abstracted from progress,
> Long and short abstracted from contrast,
> High and low abstracted from depth,
> Song and speech abstracted from melody,
> After and before abstracted from sequence.

> So it is that existence and non-existence give birth the one to (the idea of) the other (James Legge); that difficulty and ease produce the one (the idea of) the other; that length and shortness fashion out the one the figure of the other; that (the ideas of) height and lowness arise from the contrast of the one with the other; that the musical notes and tones become harmonious through the relation of one with another; and that being before and behind give the idea of one following another.
> Accordingly (back to authors), the division between truth and false comes from practice: there is truth, because the outcome is desired, and vice versa to false. Furthermore, without false, where comes the truth? And without truth, where comes the false? Human has been unintentionally comparing, weighing, balancing and trading off between favorable and unfavorable, desired and undesired, based on his prompt subjective and objective situations, hence comes the abstraction (distinction).
> Therefore, **the distinction (or abstraction) of truth and false is nothing more than desires that are subject to constant change.**

Conclusion:

**Truth is born from false and false from truth. They are exactly the measurement of men's practice. This measurement is by no means isolated from practical situations.**
- Discover the truth through practice, and again through practice verify and develop the truth. Start from perceptual knowledge and actively develop it into rational knowledge; then start from rational knowledge and actively guide revolutionary practice to change both the subjective and the objective world. Practice, knowledge, again practice, and again knowledge. This form repeats itself in endless cycles, and with each cycle the content of practice and knowledge rises to a higher level.

Accordingly, we need to represent learning in integral form:

$$Learning = \int d(perception) \cdot d(logic) = \int d(perception(object)) \cdot d(logic(object))$$

in time, space and situation domains, or in compound form:

$$Learning = \int d(perception \cdot logic)$$

since perception and logic are interchangeable and inter-transformable (they melt each other), or they are decomposable in the same manner as learning.

### 7. As a Part of Learning, Logic is Dynamic

Logic depends on perception which is subject to dynamic change with environment, i.e., **the truth value swings**.

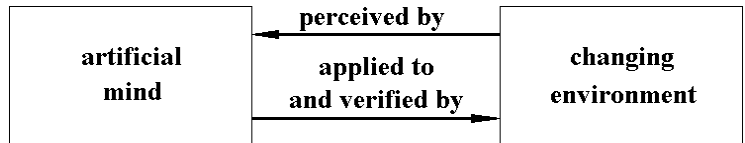

For the logic (Liu, [2]):

"I'll visit him if it doesn't rain and he is in."

To avoid being trapped in an instant case that "the clouds promise impending rain" or "there is no answer at the moment I ring the door", we need to learn the long-term trend like "does it rain whole day" or "does he keep his promise", to make a feasible plan (wait or return).

Conclusion:

$$Learning = \int d(perception \cdot logic)$$

in all its time, space and situation domains. **Whenever we persist in some instant or partial look, we loose the whole. Furthermore, whenever we are satisfied with one-sided view, we also loose it.**

### 8. As a Part of Learning, Logic is Multilateral

Human is normally too confident of himself, no matter how partial or monolateral he is, so that he always misleads himself.

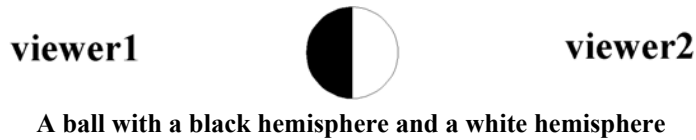

**A ball with a black hemisphere and a white hemisphere**

He may then solve such contradiction by standing on some other perspective. However, **he doesn't succeed until he reaches the opposite side**. This is where our yin-yang philosophy starts.

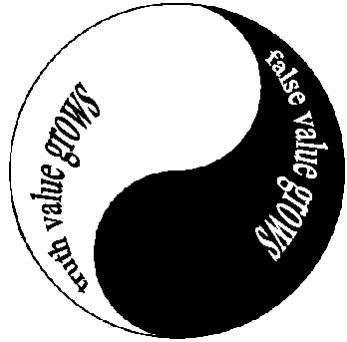 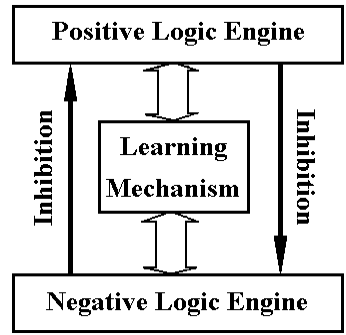

**Two kinds of partial reasoning are alternatively fired**

The incompleteness of human mind (by no means can we assert the perfectness of human being) indicates that **human always reasons in partial mode, i.e., positively partial at one time and negatively partial at another**, as indicated by taiji figure. There are plenty of reasons:
- As a holograph of the universe, human behaves in rhythmed way, e.g.: positive mode, negative mode, positive, negative… and so on, with each mode complementing and inhibiting the opposite one.
- Just because human sometimes stands on positive perspective and then the other, the truth value (it is, not it is) is in constant change, as shown in dynamic state to us.
- There remains a learning procedure hidden in the above bilateral logic: human has to balance the bilateral reasoning to adapt to his present or long-term situations, or meet his needs.
- It is from this inter-complementary and inter-inhibitory contradiction: the bilateral model, that multilateral system is generated, according to Chinese yin-yang philosophy or I-ching (in Chinese: Yijing, also known as the Book of Changes).

We can never base intelligence on the individual behavior, and this is the reason why we need group or society to exchange our views. This is also the underlying essence of multiagent approach that tries to simulate a society.

Based on our intensive exploration in Chinese philosophic perspective, a prototype of **logic cell is presented as an inter-complementary pair**: the positive and the negative logic engines represent positive partial reasoning and negative partial reasoning respectively, and a learning mechanism tries to assemble them to figure out the trend, based on its present and long-term situations, and make choices among various possibilities, e.g., to try a plan among a couple of possible plans. Suppose that each logic engine is decomposable in the same manner, which clones the entire system.

It starts with default logic—although believed absolutely valid, it is actually partial in nature as we know that human is by no means complete as long as he thinks in logic. This partial activity is warmed up and up until some time (one day) contradiction arises. This contradiction, growing up and up in previous partial mode, gives rise to negative reasoning, which later on inhibits the original logic. This inter-complementary process continues as the loop goes on and on, during which the contradiction tends to be neutralized. However, this is only a temporary balance when the two engines reach an agreement. New contradiction comes with the constant change in environment or situations outside. There still remains a chronic cycle hidden in the rhythm for long-term resolution. It is also important to note that this is a genetic proliferation in both time and space complexity, for it can clone all its subsystems.

More intensive study is being carried out in neutrosophy.

Conclusion:
A practical reference frame should originate from a single contradiction or yin-yang (see Daodejing or I-ching).
**Every existence is of bilateral character, or double characters, with each opposing and complementing the other to form a unity.**

$$d(object) = \frac{\partial(object)}{\partial (positive\ perspective)} d(positive\ engine) + \frac{\partial (object)}{\partial (negative\ perspective)} d(negative\ engine)$$

$$= \frac{\partial \text{ (object)}}{\partial \text{ (yang)}} d(yang) + \frac{\partial \text{ (object)}}{\partial \text{ (yin)}} d(yin) = \frac{d(object)}{d(reference\ frame)} d(reference\ frame)$$

Whenever we hold logic, we have already been standing on a default perspective. Is there universal logic? No, unless we reconceptualize it in an opposite perspective, e.g., Daoist or Buddhist view.

Conclusion:
When we hold logic, we have already believed that it is something. This belief in turn inhibits our negative consciousness that it may be something else (to some degree) or it can be another thing simultaneously (to some extent). We are in this way trapped. So:
> "It is never too old (for a machine) to learn."

## 9. Concluding Remarks: Illusion and Creativity

It has long been illustrated in Daodejing (Wang Bi, Guo Xiang) that whenever we capture dao as the natural law, universal method or logic, etc., what we capture can never be it. Therefore:
> **Although we can learn logic, we can never capture it.**
> **Whenever we do, what we capture is merely a distortion.**

A famous poem from "Topic to Xilin Wall" by Su Shi, a great poet in the Chinese Song Dynasty:
> A great mountain by vertical and horizontal view,
> Far, near, high, low, and each not same.
> I can't see the true face of Lushan,
> Because I am just in there.

We cannot criticize ourselves just because we habitually and absolutely believe in ourselves. Accordingly, we cannot keep a critical mind to our so called truth just because we habitually and absolutely believe in so called truth. Is there really some kind of (absolute) truth on earth?

> One time in Tang dynasty China, the Fifth Patriarch of Buddhism announced to his disciples that everyone write a verse to show his insight of the Buddhist wisdom.
> At this, the most eligible one presented on the wall the verse:
>> Our body be a Bodhi tree,
>> Our mind a mirror bright,
>> Clean and polish frequently,
>> Let no dust alight.
>
> Just as a choreman in the mill of the temple, Huineng answered it with his own:
>> There is no Bodhi tree,
>> Nor stand of a mirror bright,
>> Since all is void,
>> Where can the dust alight?

**Whenever we hold the belief "it is …", we are loosing our creativity. Whenever we hold that "it is not …", we are also loosing our creativity.** Our genuine intelligence requires that we completely free our mind - neither stick to any extremity nor to "no sticking to any assumption or belief" (Liu [3]).

As we mentioned previously, whenever there is truth, there is also false that is born from/by truth—this abstraction (distinction) is fatal to our creativity.

Meanwhile, our creativity is nothing similar with things created (i.e., in the sense if truth and false). It must be void in form (no definite form), something like dao, since whenever we hold it, it is not our creativity. Nor does it mean to destroy everything (there is nothing to destroy nor such action, if there was, it is no longer void).

> **There is nothing to destroy, nor anything to create.**
> **If there is, it is rather our illusion than our creativity.**

Because everything believed existing, true or false, is nothing more than our mental creation, there is no need to pursuit these illusions, as illustrated in the Heart Sutra [3]:
> When Bodhisattva Avalokiteshvara was practicing the profound Prajna Paramita, he illuminated the Five Skandhas and saw that they are all empty, and he crossed beyond all suffering and difficulty.

Shariputra, form does not differ from emptiness; emptiness does not differ from form. Form itself is emptiness; emptiness itself is form. So too are feeling, cognition, formation, and consciousness.

Shariputra, all Dharmas are empty of characteristics. They are not produced, not destroyed, not defiled, not pure; and they neither increase nor diminish. Therefore, in emptiness there is no form, feeling, cognition, formation, or consciousness; no eyes, ears, nose, tongue, body, or mind; no sights, sounds, smells, tastes, objects of touch, or Dharmas; no field of the eyes up to and including no field of mind consciousness; and no ignorance or ending of ignorance, up to and including no old age and death or ending of old age and death. There is no suffering, no accumulating, no extinction, and no Way, and no understanding and no attaining.

Because nothing is attained, the Bodhisattva through reliance on Prajna Paramita is unimpeded in his mind. Because there is no impediment, he is not afraid, and he leaves distorted dream-thinking far behind. Ultimately Nirvana! All Buddhas of the three periods of time attain Anuttara-samyak-sambodhi through reliance on Prajna Paramita. Therefore know that Prajna Paramita is a Great Spiritual Mantra, a Great Bright Mantra, a Supreme Mantra, an Unequalled Mantra. It can remove all suffering; it is genuine and not false. That is why the Mantra of Prajna Paramita was spoken. Recite it like this:

Gaté Gaté Paragaté Parasamgaté
Bodhi Svaha!

Conclusion:

**Everyone can extricate himself out of the maze of illusion**, said Sakyamuni and all the Buddhas, Bodhisattvas around the universe, their number is as many as that of the sands in the Ganges (Limitless Life Sutra, Chin Kung).